\documentclass[12pt]{article}
\input epsf      
\usepackage{graphicx}% Include figure files
\usepackage{amssymb}
\usepackage{amsfonts}
\usepackage{amsthm}
\usepackage{amsfonts}

\newcommand \boundary {{\partial}}

\newcommand{\integers}{{\mathbb Z}}

\newcommand \no {\noindent}
\newcommand \ra {\rightarrow}

\newcommand{\ba}[1]{\begin{array}{#1}}
\newcommand{\ea}{\end{array}}

\newcommand{\be}{\begin{equation}}
\newcommand{\ee}{\end{equation}}
\newcommand{\bea}{\begin{eqnarray}}
\newcommand{\eea}{\end{eqnarray}}
\newcommand{\beann}{\begin{eqnarray*}}
\newcommand{\eeann}{\end{eqnarray*}}

\newcommand{\weight}{{W}}

\def\reff#1{(\ref{#1})}
\begin{document}

\bibstyle{ams}

\title{Simulating self-avoiding walks in bounded domains}

\author{Tom Kennedy
\\Department of Mathematics
\\University of Arizona
\\Tucson, AZ 85721
\\ email: tgk@math.arizona.edu
}

\maketitle 

\begin{abstract}
Let $D$ be a domain in the plane containing the origin. 
We are interested in the ensemble of self-avoiding walks (SAW's) in
$D$ which start at the origin and end on the boundary of the domain. 
We introduce an ensemble of SAW's that we expect to
have the same scaling limit. The advantage of our ensemble is 
that it can be simulated using the pivot algorithm. 
Our ensemble makes it possible to accurately study SLE
predictions for the SAW in bounded simply connected domains. 
One such prediction is the distribution along the boundary of the 
endpoint of the SAW.
We use the pivot algorithm to simulate our ensemble and study this density.
In particular the lattice effects in this density that persist in the 
scaling limit are seen to be given by a purely local function. 

\end{abstract}

\section{Introduction}
\label{intro_sect}

In two dimensions there has been a lot of interest in the 
self-avoiding walk (SAW) in simply connected domains $D$ because of its
conjectured relationship with SLE. One is interested in two cases: 
the radial case in which the SAW starts at a point inside the domain and 
ends on the boundary and the chordal case in which it starts 
and ends at boundary points. If one fixes the starting and ending 
points then Lawler, Schramm and Werner conjectured that the 
scaling limit is radial or chordal SLE$_{8/3}$ \cite{lsw_saw}. 
If one fixes the starting point but allows all 
SAW's that end anywhere on the boundary, then there are conjectures 
for the hitting density along the boundary from SLE partition functions 
\cite{lsw_saw,lawler_utah}.
Some progress towards proving the conformal invariance of the SAW was made 
by Duminil-Copin and Smirnov \cite{dc_smirnov}.

For the SAW in the full plane or a half plane there is a very 
fast Monte Carlo algorithm known as the pivot algorithm. 
(Clisby's recent implementation of the 
algorithm using binary trees has dramatically increased its speed
\cite{clisby}.) However the pivot algorithm cannot
be used for the SAW in a bounded domain. The pivoting
will almost always produce a SAW that does not end on the boundary of the
domain  or leaves the domain. And the pivot algorithm works on an 
ensemble of SAW's with a fixed number of steps. For the SAW in $D$
we must allow SAW's with all numbers of steps. 

In this paper we show how to use the pivot algorithm to simulate an 
ensemble of SAW's that should have the same scaling limit as the 
ensemble of SAW's in $D$ that end on the boundary. The key idea
is to consider the ensemble of all SAW's starting at the origin 
such that there is a positive constant $\lambda>0$ such that the SAW stays
in the dilated domain $\lambda D$ and ends on the boundary of $\lambda D$. 
We refer to this ensemble as the ``dilation ensemble.'' 
The SAW's in the dilation ensemble have an arbitrary number of steps. 
Nonetheless we will argue that in the scaling limit this ensemble 
can be realized using the ensemble of walks with a fixed number of 
steps that satisfy the constraint above. 
This allows us to use the pivot algorithm to simulate the 
scaling limit of the ensemble. 

If one allows the SAW to end anywhere on the boundary of the domain, 
then in the scaling limit the distribution of the endpoint along 
the boundary gives a probability measure on the boundary analogous
to harmonic measure for the ordinary random walk. 
Assuming the boundary is smooth, the 
distribution should be absolutely continuous with respect to arc 
length along the boundary. 
Lawler, Schramm and Werner have given an explicit conjecture 
for this boundary density \cite{lsw_saw} using SLE partition functions. 
The simplest description of their conjecture
is that the boundary density for the SAW is proportional 
to the density for harmonic measure raised to the $5/8$ power. 
In \cite{lsw_saw} they 
considered domains in which the sides are parallel to a lattice 
direction. In the case of a horizontal strip with the SAW starting on one
boundary of the strip and ending on the other boundary, the
conjecture for the boundary density was tested in \cite{dyhr_et_al},
and good agreement was found. 
The conjecture was discussed further by Lawler in \cite{lawler_utah}
where it was stated for general simply connected domains 
``after taking care of the local lattice effects.''

In \cite{kennedy_lawler} it was conjectured that the local lattice
effect at a point $z$ on the boundary only depends on the 
angle of the tangent to the boundary at $z$.
There are a variety of ways to interpret what it means for the SAW 
which is on a lattice to end on the boundary of $D$ which typically
does not pass through lattice sites. 
The precise nature of the lattice effects depends on which 
interpretation is used.
Explicit conjectures for these local lattice effects were given in 
\cite{kennedy_lawler} for two particular interpretations of 
ending on the boundary. In one interpretation
one considers all SAW's $\omega$ with $\omega(i) \in D$ for $i<|\omega|$ and 
$\omega(|\omega|) \notin D$, where $|\omega|$ is the number of steps in the 
SAW. While this is one of the most natural interpretations,  
there is no good algorithm for simulating this ensemble. 
The other ensemble considered in \cite{kennedy_lawler}
is the ``cut-curve'' ensemble.
It uses infinite length SAW's in the full plane and conditions
on the event that the SAW crosses the boundary of $D$ only once.
(In practice we use SAW's with a fixed length $N$ and take $N$ large
enough that the typical size of the SAW is much larger than the 
domain $D$.) SAW's in this ensemble are equivalent to the concatenation
of a SAW in $D$ from the origin to the boundary of $D$
and a SAW in the exterior of $D$, starting at the same boundary point 
and going to $\infty$.

One can use the pivot algorithm to simulate the ``cut-curve'' ensemble. 
But this has the disadvantage 
that one must deal with a double limit - one must first let the 
number of steps in the SAW go to infinity and then take the lattice 
spacing to zero. In practice this means that one must simulate 
SAW's that are much larger than the domain $D$. So such simulations 
are less efficient than the method we will introduce for the 
dilation ensemble. 

In the next section we define the dilation ensemble and show how it 
may be simulated using the pivot algorithm on the ensemble of walks 
with a fixed number of steps. In section three we give an explicit 
conjecture for the lattice effects in our dilation ensemble  and show 
how these effects may be computed by simulation. In section four we 
use our method for simulating the dilation ensemble to study the SAW
in several different domains. We compare the boundary densities we find
with the SLE partition function predictions for these densities.

\section{The dilation ensemble}
\label{dilation_ensemble_sect}

We first consider the radial case. 
We assume that our simply connected domain $D$ contains the origin
and has a boundary which is a smooth curve which we denote by $\boundary D$. 
We also assume that the domain is star-shaped with respect to 
the origin, i.e., any ray from the origin intersects $\boundary D$ in only 
one point. (It would be interesting to generalize our approach 
to domains that do not have this property.) 
We are interested in the ensemble of SAW's that start 
at the origin and end on $\boundary D$.
Note that we are not fixing the point on the boundary where the SAW must
end; the endpoint of the SAW is a random point on the boundary.

The dilation ensemble can be described as all 
SAW's such that the domain $D$ can be dilated so that the SAW lies in the 
dilated domain and ends on the boundary of the dilated domain. 
We let $\delta>0$ denote the lattice spacing. 
We start with the ensemble of all SAW's in the full plane on the 
lattice $\delta \integers^2$ 
that start at the origin with any finite number of steps.
For a SAW $\omega$ we let $|\omega|$ denote the number of steps in 
the walk. The walk $\omega$ is given the weight $\mu^{-|\omega|}$ 
where $\mu$ is the connective constant. 
(It is given by $\mu = \lim_{N \ra \infty} c_N^{1/N}$ where 
$c_N$ denotes the number of SAW's in the full plane with $N$ steps that 
start at $0$. The existence of this limit follows from the 
subadditivity of $\ln(c_N)$ \cite{madras_slade}.)
We will impose two constraints on the ensemble of all finite length SAW's. 
We let $\lambda(\omega)>0$ be the dilation
such that the endpoint of $\omega$ is on the dilated 
curve $\lambda(\omega) \boundary D$. 
The partition function for the ensemble of all finite length SAW's 
starting at the origin is infinite, so  
we add the constraint that 
$\lambda(\omega)$ lies in $[\lambda_1,\lambda_2]$ where  
$\lambda_1 < \lambda_2$ are positive constants.  
We expect that this makes the partition function finite, although 
we cannot prove this. 
The walk $\omega$ ends on the curve $\lambda(\omega) \boundary D$, but it 
need not stay inside the region $\lambda(\omega) D$.
Our second constraint is that the walk $\omega$ does stay strictly inside
this region except for the endpoint.
We let $1_D(\omega)$ denote the indicator function that is $1$ if
$\omega$ stays strictly inside this region except for its endpoint and 
$0$ if it does not. 
So our ensemble corresponds to the partition function 
\bea
Z= \sum_{\omega} \mu^{-|\omega|} 1_D(\omega) \, 
1(\lambda_1 \le \lambda(\omega) \le \lambda_2) 
\, \weight(\omega)
\label{uncoiled_z}
\eea
We have included a weighting factor $\weight(\omega)$ which we will 
define later to make our ensemble correspond to the ensemble of 
SAW's in $D$ that end on the boundary.
Associated with this partition function is a probability measure
on simple curves in $D$ that end on the boundary.
Let $\psi(\gamma)$ be a function on simple curves $\gamma$  in $D$
that start at the origin and end on $\boundary D$. 
If $1_D(\omega)=1$, then ${\omega \over \lambda(\omega)}$ is such a curve
and so $\psi({\omega \over \lambda(\omega)})$ is defined. 
Define 
\bea
Z(\psi)= \sum_{\omega} \mu^{-|\omega|} 1_D(\omega) \, 
1(\lambda_1 \le \lambda(\omega) \le \lambda_2) \, \weight(\omega) \,
\psi({\omega \over \lambda(\omega)})
\label{z_psi}
\eea
The expected value of $\psi$ is defined to be $Z(\psi)/Z$. 
(Both $Z$ and $Z(\psi)$ depend on $\lambda_1$ and $\lambda_2$, but they
will be fixed throughout our discussion.) 

We now consider how this ensemble corresponds to the ensemble of 
SAW's in $D$ that end on $\boundary D$. 
For a large integer $N$, let $d \lambda= (\lambda_2-\lambda_1)/N$.
We think of the region between the curves $\lambda_1 \boundary D$ 
and $\lambda_2 \boundary D$ 
as the union of the regions between the curves 
$(\lambda_1 + (k-1) d \lambda)\boundary D$ and
$(\lambda_1 + k d \lambda)\boundary D$ where $k=1,2,\cdots,N$. 
In the scaling limit, the ensembles corresponding to these different
regions are related by just a dilation. 
So they will all give the same expected value to $\psi$. 
The region between $(\lambda_1 + (k-1) d \lambda)\boundary D$ and
$(\lambda_1 + k d \lambda)\boundary D$ corresponds to thickening the boundary 
of $(\lambda_1 + k d \lambda)\boundary D$. So we can think of the
dilation ensemble as interpreting the constraint that the SAW stays
in $D$ and ends on $\boundary D$ as thickening the boundary. 

It is important to observe that the way that our ensemble thickens the 
boundary is not uniform along $\boundary D$. 
The natural way to thicken the curve would be to take the thickness 
in the direction perpendicular to the curve to be constant along the curve. 
Let $\rho(z)$ denote the density of the endpoint of the SAW in 
$D$ which ends on the boundary.  In the dilation ensemble if we 
take the weight $\weight(\omega)$ to just be constant, then the 
corresponding boundary density will be proportional to $\rho(z)$ times
the thickness at $z$. To correct for this we take 
$\weight(\omega)$ to be proportional to the inverse of this thickness. 

Recall that we assume that our domain is star-shaped with respect 
to the origin, i.e., any ray from the origin intersects $\boundary D$ in 
exactly one point. So we can parametrize $\boundary D$ by the polar 
angle $\theta$ with respect to the origin.  We let $D(\theta)$
be the distance from the origin to the curve in the 
direction $\theta$. 
Along the ray at polar angle $\theta$ the thickness is 
$(\lambda+d\lambda) D(\theta) - \lambda D(\theta) = D(\theta) d \lambda$. 
But this segment is not perpendicular to the tangent line. 
Let $\alpha(\theta)$ be the angle of a line perpendicular 
to $\boundary D$ at our point.
Then the thickness of our shell perpendicular to the tangent line is
$D(\theta) \cos(\theta-\alpha(\theta))$. 
Thus we define 
\bea
\weight(\omega) =  [D(\theta(\omega)) 
\cos(\theta(\omega)-\alpha(\theta(\omega))) ]^{-1}
\label{weight_def}
\eea
where $\theta(\omega)$ is the polar angle of the endpoint of $\omega$.
Note that the weight only depends on the polar angle of the 
endpoint of $\omega$. 

The dilation ensemble \reff{uncoiled_z} includes SAW's of all lengths. 
To simulate it using the pivot algorithm we must relate it to the 
ensemble of SAW's with a fixed length in which each SAW has the 
same weight. Our method for doing this is closely related to the method in 
\cite{kennedy} for relating the ensemble of SAW's in the half-plane with 
a fixed number of steps to radial SLE. 
We decompose the sum in $Z(\psi)$, 
as defined in \reff{z_psi}, according to the length of $\omega$. 
\beann
Z(\psi) = \sum_n \mu^{-n} \, \sum_{\omega:|\omega|=n} 
1_D(\omega) \, 
1(\lambda_1 \le \lambda(\omega) \le \lambda_2) 
\, \weight(\omega) \, \psi({\omega \over \lambda(\omega)}) 
\eeann
Let $c_n$ be number of SAW in the full plane starting at $0$ 
with $n$ steps. We have
\beann
Z(\psi) = \sum_n c_n \mu^{-n} \, {1 \over c_n}
\sum_{\omega: |\omega|=n} \, 1_D(\omega) \, 
1(\lambda_1 \le \lambda(\omega) \le \lambda_2) \, 
\weight(\omega) \, \psi({\omega \over \lambda(\omega)}) 
\eeann
The constraint $\lambda(\omega) \ge \lambda_1$ implies that $\omega$ must
have at least $\lambda_1/\delta$ steps. So as the lattice spacing goes 
to zero, the first $n$ for which the summand in the sum on $n$ is nonzero
goes to infinity. 
Since $c_n$ is asymptotic to $\mu^{n} n^{\gamma-1}$, 
we replace $c_n \mu^{-n}$ by $n^{\gamma-1}$.
Let $P_n$ be the uniform probability measure on all SAW's starting at the 
origin with $n$ steps, and let $E_n$ be the associated expected value. 
Then we can write the above as 
\beann
Z(\psi) = \sum_n n^{\gamma-1} \, E_n[1_D(\omega) \, 
1(\lambda_1 \le \lambda(\omega) \le \lambda_2) \, 
\weight(\omega) \, \psi({\omega \over \lambda(\omega)}) ]
\eeann

As noted before, the constraint $\lambda(\omega) \ge \lambda_1$ 
implies that the sum on $n$ is only over large values. So we should be 
able to approximate $E_n$ by its scaling limit. The constraint
$1_D(\omega)$ is a bit tricky in the scaling limit, so we proceed
as follows. Fix a large positive integer $N$. If $n$ is also large
and $\omega$ and $\gamma$ are drawn from $P_n$ and $P_N$, respectively, then  
the distributions of $\delta^{-1} n^{-\nu} \omega$ and 
$\delta^{-1} N^{-\nu} \gamma$ are approximately the same. 
So we will replace $E_n$ by $E_N$ by replacing 
$\delta^{-1} n^{-\nu} \omega$ by $\delta^{-1} N^{-\nu} \gamma$,
i.e., we replace $\omega$ by $n^\nu N^{-\nu} \gamma$.
So $\lambda(\omega)$ becomes 
$\lambda(n^\nu N^{-\nu} \gamma) = n^\nu N^{-\nu} \lambda(\gamma)$.
Our weight $\weight(\omega)$ only depends on the polar angle of 
the endpoint of $\omega$, 
so we can replace $\weight(\omega)$ by $\weight(\gamma)$. 
And we can replace $\psi({\omega \over \lambda(\omega)})$ by 
$\psi({\gamma \over \lambda(\gamma)})$.

The constraint $1_D(\omega)$ is more subtle. 
The probability that an $n$-step SAW stays on one side of 
a half plane is conjectured to go to zero as $n^{-\rho}$ as $n \ra \infty$
with $\rho=25/64$ \cite{lsw_saw}.
So we expect that the probability that  $1_D(\omega)=1$ also goes
to zero as $n^{-\rho}$.
So we approximate $n^\rho 1_D(\omega)$ with $N^\rho 1_D(\gamma)$, i.e., 
we replace $1_D(\omega)$ by $N^\rho n^{-\rho} 1_D(\gamma)$.
We now have
\beann
Z(\psi) \approx N^\rho \sum_n n^{\gamma-1-\rho} \, E_N [ 1_D(\gamma) 
\,  1(\lambda_1 \le n^\nu N^{-\nu} \lambda(\gamma) \le \lambda_2) 
\, \weight(\gamma) \, \psi({\gamma \over \lambda(\gamma)}) ]
\eeann

The $n$ dependent part of this is 
\beann
\sum_n n^{\gamma-1-\rho} \, \,
 1(\lambda_1 \le n^\nu N^{-\nu} \lambda(\gamma) \le \lambda_2)
\eeann
Since the constraint restricts the sum to large values of $n$, 
$n^{\gamma-1-\rho}$ is slowly varying and so we can think of this as 
a Riemann sum approximation to 
\beann
\int_0^\infty x^{\gamma-1-\rho} \, 
1(\lambda_1 \le x^\nu N^{-\nu} \lambda(\gamma) \le \lambda_2)
\, dx = c \, N^{\gamma-\rho} [\lambda(\gamma)]^{(\rho-\gamma)/\nu}
\eeann
where the constant $c$ depends on $\lambda_1$ and $\lambda_2$, but 
nothing else. 
Thus
\beann
Z(\psi) \approx c N^{\gamma-\rho} \,
E_N [ \lambda(\gamma)^{(\rho-\gamma)/\nu} \, 1_D(\gamma) \, 
\weight(\gamma) \, \psi({\gamma \over \lambda(\gamma)}) ]
\eeann
By taking $\psi=1$, this result also gives an expression for $Z$. 

The above derivation was for the radial case in which the starting point 
of the SAW is in the interior of $D$. 
If we take a domain which has the origin on its boundary and consider
SAW's which start at the origin, stay in $D$ and end on its
boundary, then we can repeat the above derivation. 
The one important difference is that the number of walks $c_n$ 
should now grow like $\mu^n n^{\gamma-1-\rho}$, rather than 
$\mu^n n^{\gamma-1}$. So the appropriate
power for $\lambda(\gamma)$ becomes $p= {2\rho-\gamma \over \nu }$.

We conjecture that when we take the scaling limit 
our approximations become exact. More precisely, we make the 
following conjecture.  

\medskip
\no {\bf \large Conjecture :} 
{\it 
Let $D$ be a simply connected domain which contains the origin
(radial case) or has the origin on its boundary (chordal case), and 
which is star shaped with respect to the origin. Let $\psi$ be a function
on simple curves in $D$ which start at the origin and end on the boundary.  
Let $Z(\psi)/Z$ be the expected value of $\psi$ in the dilation
ensemble of SAW's on a lattice of spacing $\delta$ as defined by \reff{z_psi}.
Let $E_N$ be the uniform probability measure on $N$-step SAW's 
$\gamma$ starting at the origin in the half plane (chordal case) or 
in the full plane (radial case).
Let $\lambda(\gamma)>0$ be such that $\gamma$ 
ends on $\lambda(\gamma) \boundary D$. Let $1_D(\gamma)$ be 
the indicator function of the event that $\gamma$ is 
inside the domain $\lambda(\gamma) D$ except for its endpoint(s). 
Define the  weight $\weight(\gamma)$ by \reff{weight_def}. Then 
\bea
\lim_{\delta \ra 0} {Z(\psi) \over Z} =
\lim_{N \ra \infty} { E_N[ \lambda(\gamma)^p \, 1_D(\gamma) \, 
\weight(\gamma) \, \psi({\gamma \over \lambda(\gamma)}) ]
\over 
E_N [ \lambda(\gamma)^p \, 1_D(\gamma) \, \weight(\gamma)]}
\label{conjecture}
\eea
\beann
p &=& {\rho-\gamma \over \nu} =-61/48 \quad (radial \, case),  \\
p &=& {2\rho-\gamma \over \nu } =-3/4 \quad (chordal \, case).
\label{pvalues}
\eeann
}

\section{Lattice effects in boundary densities}
\label{lattice_effect_sect}

We expect that in the scaling limit 
the endpoint of the SAW on the boundary of the domain $D$ 
will have a distribution that is absolutely continuous with 
respect to arc length along the boundary.
We denote this boundary density by $\rho$.
If we consider ordinary random walks instead of 
self-avoiding walks, then in the scaling limit this boundary 
distribution would be harmonic measure. 

Let $f$ be a conformal map on $D$. Then it is expected that the boundary
density transforms as 
\bea
\rho_{D}(z) = c |f^\prime(z)|^{5/8} \, \rho_{f(D)}(f(z)), \quad z \in 
\boundary D 
\label{conf_covar} 
\eea
except for a local lattice effect that persists in the scaling limit. 
If we take $g$ to be the conformal map of $D$ onto the unit disc 
that fixes the origin, then the above implies that 
$\rho_{D}(z)$ is proportional to $|g^\prime(z)|^{5/8}$.
Note that the conformal invariance of harmonic measure implies
that the boundary density for the ordinary random walk 
transforms in an analogous way
except that the power of $5/8$ is replaced by $1$. 
Consequently \reff{conf_covar} implies that the boundary density for 
the SAW is proportional to the boundary density for the ordinary 
random walk raised to the $5/8$ power.

We now turn to the computation of the lattice effect correction to 
the boundary density.
This computation for the dilation ensemble is analogous to 
what was done in \cite{kennedy_lawler}.
The constraint that the SAW $\omega$
stays inside the dilated curve $\lambda(w) \boundary D$ 
has both a macroscopic and microscopic nature. 
The conjecture \reff{conf_covar} is the result of the 
macroscopic effect. 
Near the endpoint of the walk there is an additional microscopic effect.
Consider a SAW that ends at $z$ and consider the tangent line to the 
curve $\boundary D$ at $z$. The constraint that the SAW stays inside 
$\lambda(\omega) D$ implies that near the endpoint 
the SAW must stay on one side of this line. 
This will produce a factor $l(\theta)$ that depends on the angle of 
the tangent line with respect to the lattice. 

For an angle $\theta$, let $L$ be the line with 
polar angle $\theta$ passing through the origin. 
We consider walks with $N$ steps starting at the origin.
Let $c_N$ be the number of such walks, and let
$b_N(\theta)$ be the number of such walks that stay on one side of the line.
So $b_N(\theta)/c_N$ is the probability that an 
$N$ step SAW stays on one side of the line. We expect 
that this probability goes to zero as $N^{-\rho}$ as $N \ra \infty$, 
and we conjecture that the lattice effect function is given by 
\bea
l(\theta)= \lim_{N \ra \infty} {b_N(\theta) \over c_N} N^{\rho}
\eea
Since $\theta$ and $\theta+180$ give the same line, 
$l(\theta)$ has period $180$ degrees. 
The function will have more symmetries depending on the type of lattice. 
For example, for the square lattice $l(\theta)$
has period $90$ degrees and  $l(\theta)=l(90-\theta)$. 
If we use the weight $\weight(\omega)$ in \reff{uncoiled_z} and 
\reff{z_psi}, then the boundary density for 
the dilation ensemble will be $\rho_D(z) l(\theta(\omega))$ where
$\theta(\omega)$ the polar angle of the endpoint of $\omega$.
To remove this lattice effect from our dilation ensemble we replace
the weight $\weight(\omega)$ given by \reff{weight_def} by 
\bea
\hat{\weight}(\omega)={\weight(\omega) \over l(\theta(\omega))}
\label{weight_hat}
\eea
In our conjecture \reff{conjecture} the weight $\weight(\omega)$ is 
used in both sides of this equation. The same derivation shows that this
equality should also hold if we use $\hat{\weight}(\omega)$ in 
both sides.

We take the constraint that $\omega$ stays inside $\lambda(\omega) D$
to mean that it stays strictly inside the curve except for the endpoint. 
The other convention would be to allow $\omega$ to have sites that lie
on $\lambda(\omega) \boundary D$ in addition to its endpoint. 
If $\boundary D$ has 
flat sections, then there can be a big difference between these two conventions
for angles $\theta$ such that a line through the origin at angle $\theta$
passes through some lattice sites. This will be the case for the 
equilateral triangle that we consider in our simulations.

\begin{figure}[tbh]
\includegraphics{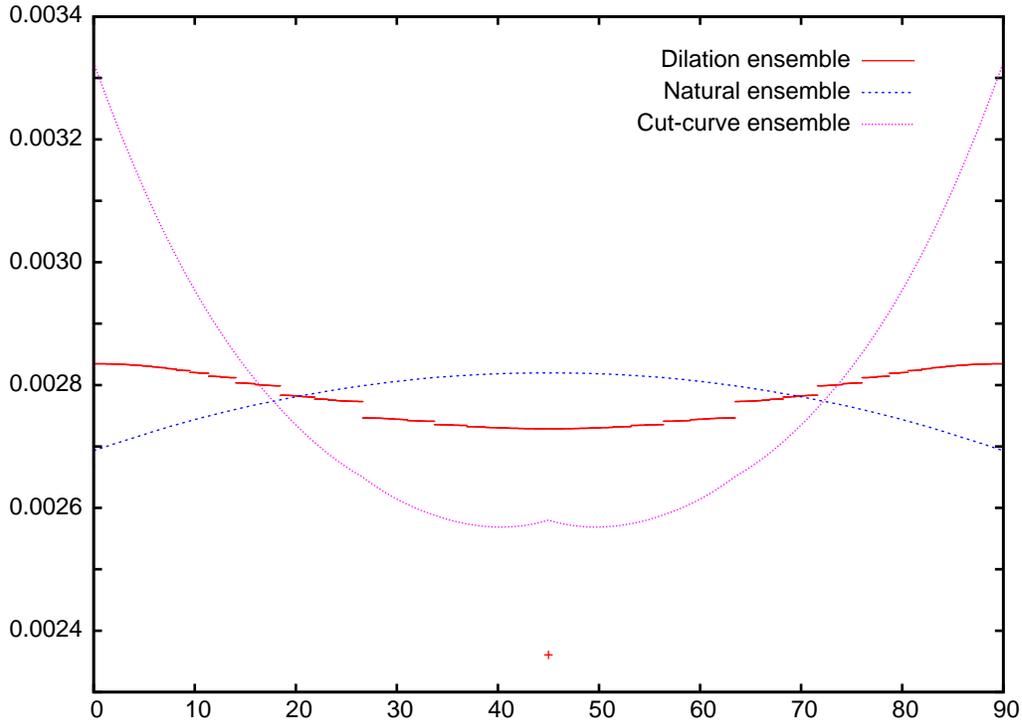}
\caption{\leftskip=25 pt \rightskip= 25 pt 
The lattice effect function $l(\theta)$ for three different ensembles
on the square lattice.}
\label{ltheta}
\end{figure}

In figure \ref{ltheta} we plot the lattice effect function $l(\theta)$
for the dilation ensemble and for two other ensembles, all on the 
square lattice. One is the 
cut-curve ensemble studied in \cite{kennedy_lawler}. The other ensemble 
is the ``natural'' ensemble in which we take all SAW's such that 
all the sites on the walk are in $D$ except for the endpoint which is 
outside of $D$. In other words, the last bond of the walk crosses the 
boundary of $D$ and is the only bond which does so. 
All three functions have been normalized so that the total 
area under each curve is $1$. 
It is important to note that the vertical scale of the plot does not 
include $0$. These functions are actually relatively flat. 
Note that the lattice effect function for the dilation ensemble 
studied in this paper is considerably flatter than the function for
the cut-curve ensemble that was simulated in \cite{kennedy_lawler}.

The function $l(\theta)$ is continuous for the cut-curve and natural 
ensembles, but for the dilation ensemble it is not. 
For the square lattice it is discontinuous
at $\theta$ such that $\tan(\theta)$ is rational. 
The biggest discontinuity is at $\theta=0$, but it cannot be seen in the 
figure since $l(0)=0.001516$ which is well below the region plotted. 
The second biggest discontinuity at $\theta=45$ is seen in the figure 
as an isolated point at $l(45)=0.002361$. 
This discontinuity is atypical in that 
$l(45-)=l(45+)$. Other discontinuities with $l(\theta-) \neq l(\theta+)$ 
can be seen in the figure. 
To see why $l(\theta)$ is discontinuous, consider the case of $\theta=0$.
Recall that we make the convention that the walk must stay strictly 
inside the domain except for the endpoint. So when $\theta$ is exactly $0$, 
$b_N(\theta)$ counts walks that start at the origin and then stay 
strictly on one side of the horizontal axis. So the walk cannot visit 
any site on the horizontal axis. Now consider a $\theta$ that is slightly 
greater than zero and consider a walk which stays above this line. 
The sites on the negative horizontal axis are now above the line and so 
the walk can visit these sites. Thus $b_N(0)$ is significantly smaller
than $b_N(0^+)$.
If we change the convention that the SAW must stay strictly inside 
the domain to allow SAW's that visit sites on the boundary, then the 
values at the angles where $l(\theta)$ is discontinuous will change 
but it will still be discontinuous.

\section{Simulations}

We use our conjecture to simulate the dilation ensemble 
and compare the boundary density found in the simulation with 
the density given by \reff{conf_covar}. Note that our simulations are 
testing three different conjectures.  One is the conjecture 
\reff{conjecture} that says we can use the fixed length ensemble to 
simulate the ensemble of SAW's in a bounded domain. 
Another is the SLE partition function prediction for the boundary 
density \reff{conf_covar}. 
And finally there is the conjecture that the lattice 
effect correction to this density is given by $l(\theta)$, i.e., 
by using the weight \reff{weight_hat}.

The exact predictions for the boundary density from SLE partition
functions \reff{conf_covar} are given in the appendix. 
In our simulations we work with cumulative distribution functions
(cdf's) instead of densities. Extracting a density from the simulation
requires taking a numerical derivative, i.e., choosing a bin size and 
computing a histogram. Using the cdf's avoids this. 

In our conjecture \reff{conjecture} 
the probability measure $E_N$ only depends on 
whether the geometry we are studying is chordal or radial. 
So we only need to do two simulations, one for the chordal cases and 
one for the radial cases. In each simulation we sample the Markov 
chain every $100$ time steps. In the chordal case we generated 
$9 \times 10^9$ samples and in the radial case $11 \times 10^9$ samples. 
For each sample and each domain $D$ we test if $1_D(\omega)=1$.
When it does we use that sample in the computation of the cdf of 
the boundary density for that domain. 
The probability that $1_D(\omega)=1$ depends on $N$ and on the domain.
In our simulations it ranges from ${1 \over 3} \%$ to $2 \%$. 

The first domain we consider is a horizontal strip of height $1$.
For this domain there are no lattice effects since the entire 
boundary is horizontal. 
We consider both a chordal case and a radial case. In the chordal 
case the strip is $\{z : 0 < Im(z) < 1\}$. The SAW starts at the 
origin and ends on the upper boundary. 
In the radial case the strip is $\{z : -1/4 < Im(z) < 3/4 \}$. 
The SAW starts at the origin and can end on either boundary.
The conjectured density for the chordal case was tested by simulation in 
\cite{dyhr_et_al}. Here we are primarily interested in using 
these domains to test our conjecture \reff{conjecture}
that uses the fixed length ensemble for the simulation, especially 
the value of the power $p$. 

\begin{figure}[tbh]
\includegraphics{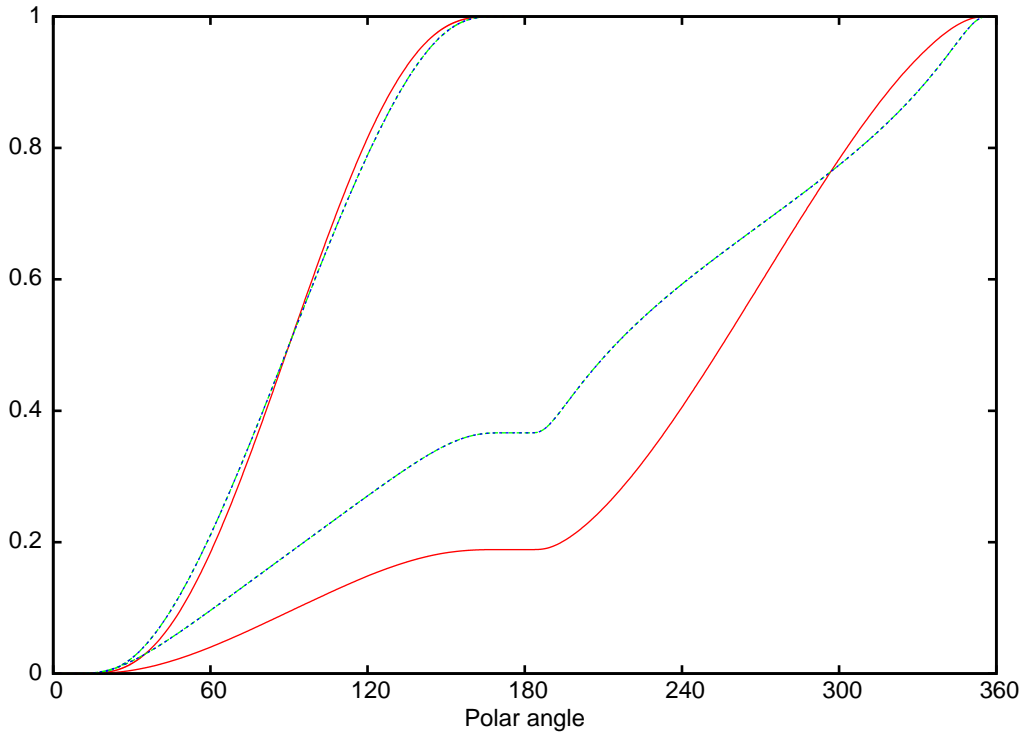}
\caption{\leftskip=25 pt \rightskip= 25 pt 
Comparison of cdf's from simulations and theory for chordal strip 
(polar angle ranges over $[0,180]$) and radial strip (polar angle ranges 
over $[0,360]$). The solid curves are simulations using $p=0$. 
The overlying dashed curves are the simulations with the correct $p$ and 
the SLE prediction.
}
\label{strip_cdf}
\end{figure}

In figure \ref{strip_cdf} we plot six curves. Two curves are 
the cdf's computed using SLE partition functions for the 
chordal and radial cases. The other four curves are simulation cdf's
for the chordal and radial cases. In each case we have two simulations, 
one with the power $p$ given by \reff{pvalues} and one with $p=0$. 
The curves from the SLE prediction and the curves from simulations 
with the correct value of $p$ are indistinguishable on this plot,
and so it appears there are only four curves in the figure. 
The maximum difference between these curves is given in table 
\ref{maxcdf}. It is on the order of $10^{-4}$ in the chordal case
and $1.5 \times 10^{-4}$ in the radial case. 
The discrepancy for the curves with $p=0$ is quite large for the 
radial strip. For the chordal strip it is smaller, but can still 
be clearly seen in the plot. The maxima of the differences for 
the $p=0$ case is also given in the table.

\begin{table}
\begin{center}
\begin{tabular}{|c|c|c|c|}
\hline
Domain & & $l(\theta)$ used ? & max of difference of cdf's \\
 & & & (in thousandths)\\
\hline
\hline
 Strip with $p=0$  & chordal &     & $ 26.90980 \pm   0.63765 $  \\
\hline
 Strip with $p=0$  & radial  &     & $  229.35605 \pm   0.64829 $ \\
\hline
 Strip             & chordal &     & $  0.09842 \pm   0.11619 $  \\
\hline
 Strip             & radial  &     & $  0.15563 \pm   0.15085 $  \\
\hline
 Triangle          & radial  & yes & $  0.05939 \pm   0.08689 $  \\
\hline
 Centered circle   & radial  & no  
                 &     $  3.16364 \pm   0.07145 $  \\
 Centered circle   & radial  & yes
           &     $  0.10997 \pm   0.07793 $  \\
\hline
 Off-center circle & radial  & no  
           &    $   2.26140 \pm   0.09184 $ \\
 Off-center circle & radial  &  yes
           &    $   0.07513 \pm   0.08061 $ \\
\hline
 Partial circle    & chordal &  no  
           &     $  1.28211 \pm   0.06962 $  \\
 Partial circle    & chordal &  yes
           &     $  0.09512 \pm   0.06158 $  \\
\hline
 Tangent circle    & chordal &  no  
           &     $  1.87370 \pm   0.10155 $  \\
 Tangent circle    & chordal &  yes
           &     $  0.18705 \pm   0.10116 $  \\
\hline
\end{tabular}
\caption{\leftskip=25 pt \rightskip= 25 pt 
For each domain the last column is the maximum of the absolute
value of the difference between two cdf's. 
One cdf is computed using the SLE partition function prediction 
for the density. The other is from the SAW simulation. 
The middle column indicates if the lattice effect $l(\theta)$ is 
corrected for.
}
\label{maxcdf}
\end{center}
\end{table}

We can estimate the power $p$ by minimizing the difference 
between the boundary density computed by simulation using the 
conjecture \reff{conjecture}
and the SLE prediction for the boundary density. We use the two 
strip geometries described above and minimize the $L^2$ norm of 
the difference between the densities. 
Our simulations to estimate the power $p$ are separate from the simulations
to compute the boundary densities. We use SAW's with $100,000$ 
steps. For the chordal geometry we generated approximately $8$ million 
samples and for the radial geometry approximately $10$ million samples.
By samples we mean SAW's for which $1_D(\omega)=1$. 
In the chordal case the conjectured exact value is 
$p= -3/4 = -0.75$  and from the simulations we find the 
minimum difference is when $p=-0.751874$, a difference of $0.25 \%$.
In the radial case the conjectured exact value is 
$p= -61/48 = -1.27083\bar{3}$ and from the simulations we find 
$-1.269917$, a difference of $0.07 \%$.  

\begin{figure}[tbh]
\includegraphics{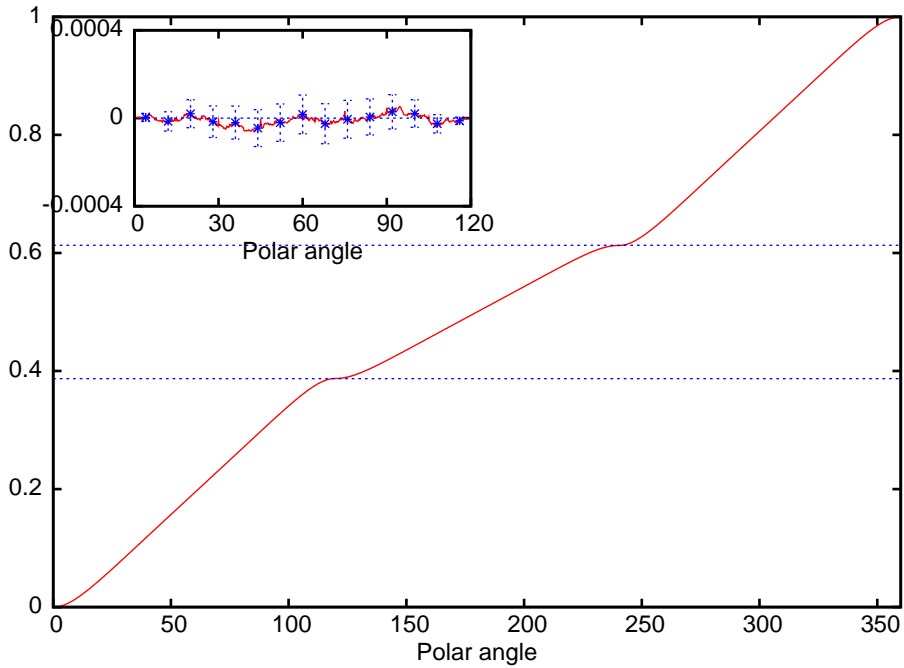}
\caption{\leftskip=25 pt \rightskip= 25 pt 
Simulation cdf for the equilateral triangle. 
The inset shows the difference of the simulation cdf and the theoretical 
cdf when the polar angle is reduced mod 120 to remove the lattice effects.
}
\label{eq_tri_cdf}
\end{figure}

Next we consider an equilateral triangle centered at the origin 
whose vertices have polar angles of $0,120$ and $240$. 
The side corresponding to polar angles in $[120,240]$ is vertical and 
the sides corresponding to $[0,120]$ and $[240,360]$ are at $30$ degrees
with respect to the lattice directions.
So two of the sides will have the same value of $l(\theta)$ 
while the third side has a different value.
This geometry gives an extreme example of the lattice effects. 
We simulate the ensemble with the weight factor $\weight(\omega)$ 
which does not correct for the lattice effects.
We find that the probabilities for hitting the 
sides corresponding to $[0,120],[120,240],[240,360]$ are 
$0.387375, 0.225173,0.387452$ respectively. 
The ratio of the smaller probability to the average of the other 
two probabilities is $0.581221$. This should be compared with the 
ratio of the two values of the lattice effect function which 
is $l(0)/l(30) =0.581281$.
In figure \ref{eq_tri_cdf} we plot the cdf from our simulation. 
The smaller probability of hitting the vertical side is clearly seen. 
The two horizontal lines are at the heights predicted by the lattice
effect function. For this geometry the only lattice effect should
be to make the probabilities of hitting the three different sides 
unequal. Given that you hit a particular edge, the distribution along 
that edge should be the same for the three edges. 
So we can remove the lattice effect by looking at the polar angle of 
the endpoint mod $120$. We compute the cdf of this random variable. 
In the inset in figure \ref{eq_tri_cdf} we show this cdf minus the SLE
partition function prediction for this cdf. This difference is quite
small, on the order of $5 \times 10^{-5}$. The maximum of the difference
is given in table \ref{maxcdf}.

\begin{figure}[tbh]
\includegraphics{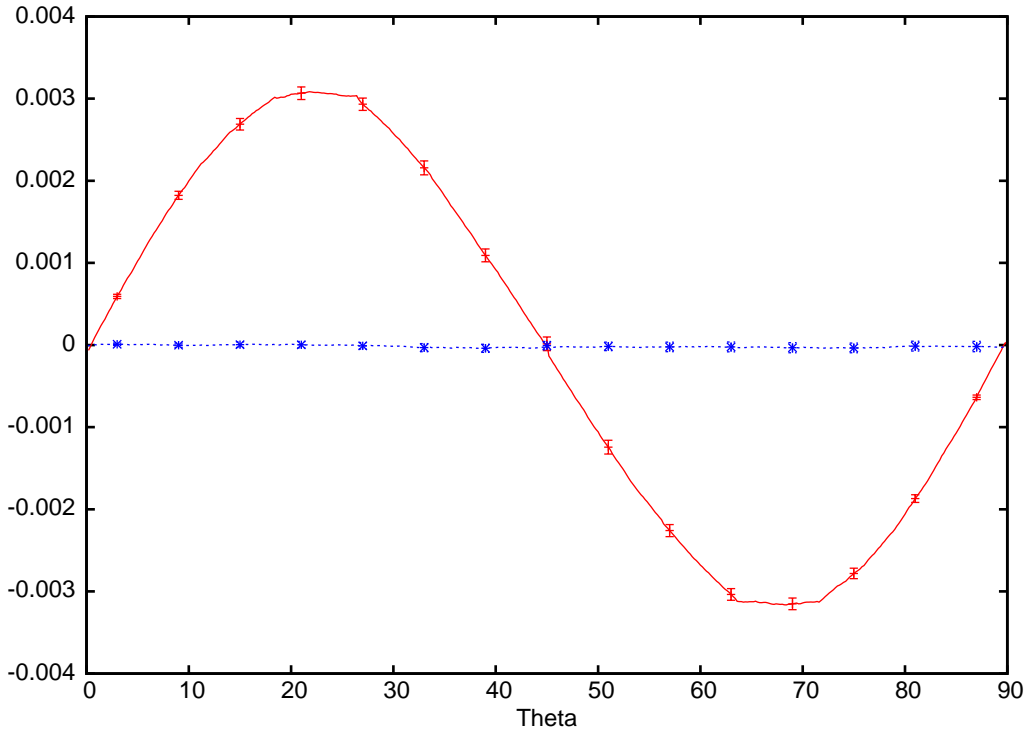}
\caption{\leftskip=25 pt \rightskip= 25 pt 
Difference of cdf's from simulations and theory for disc with origin 
at the center. 
The sine-like curve is for the simulation which does not correct for 
the lattice effects; the flat curve is for the simulation that does correct.
}
\label{circle_dif_lat}
\end{figure}

The next domain we consider is a unit disc centered at the origin. 
As always, the SAW starts at the origin. 
The SLE partition function prediction is just that the density is uniform. 
The deviation from this is entirely due to lattice effects. 
The boundary density will have a period of $90$ degrees, so we take the 
polar angle of the endpoint of the walk mod $90$. 
In figure \ref{circle_dif_lat} we plot two curves. 
Both are the cdf from the simulation minus the cdf of the uniform 
distribution. In one we use the 
weight factor $\weight(\omega)$ which does not include the lattice 
effect correction, and in the other we use the weight factor 
$\hat{\weight}(\omega)$ which does include the correction.
The sine-like curve is for the simulation without the 
lattice correction. 
The difference is small, on the order of $0.003$, but clearly not zero.
The flat curve is the difference when we do include the lattice effect 
correction in the simulation. 
The size of this difference is a test of our prediction for the lattice
effect correction.
The difference is extremely small, on the order of $10^{-4}$.
The maximum of the difference with and without the lattice correction 
is given in table \ref{maxcdf}.

Our last three simulations involve geometries for which the 
conjectured density from \reff{conf_covar} is non-trivial and 
the lattice function $l(\theta)$ enters in a non-trivial way.
In the first of these geometries the domain is a unit disc centered 
at $3/4$. We consider
SAW's that start at the origin and end on the arc of the circle.  
As before we consider the difference between the cdf 
from the simulation and the cdf predicted by SLE partition functions. 
As before we consider two cases, one using the weight $\weight(\omega)$ 
and one using the weight $\hat{\weight}(\omega)$.
We refer to this geometry as the ``off-center circle'' and 
give the  maximum of the difference of the cdf's in table \ref{maxcdf}.
As in the previous geometry the difference when we do not 
correct for the lattice effect is clearly non-zero, 
while it is zero within the errors in the 
simulation when we do correct for the lattice effects.

In the next geometry the domain is the intersection 
of the unit disc centered at $-3i/4$ and the upper half plane. 
So the domain is bounded by the real axis and the portion of the unit
circle centered at $-3i/4$ that lies above the real axis. We consider
SAW's that start at the origin and end on the arc of the circle.  
This geometry is called ``partial circle'' in the table. 
In the final geometry the domain is a unit disc centered 
at $i$. So the real axis is tangent to the disc at the origin.
We consider SAW's starting at the origin and ending on the upper 
half of the boundary. So the polar angle with respect to the origin
ranges from $45$ to $135$ degrees. 
This geometry is called ``tangent circle'' in the table.
For both of these geometries the maxima of the differences is 
given in table \ref{maxcdf}. 
They are clearly non-zero when we do not correct for the lattice effects, 
while they are zero within the errors in the 
simulation when we do correct for the lattice effects.

\section{Conclusions}

We have studied SAW's in a bounded domain $D$ which start at a fixed
point (either in the interior or on the boundary) and end 
at an unconstrained point on the boundary. 
The dilation ensemble of SAW's consists of all SAW's starting at the 
origin with the property that they can be scaled to give a curve
inside our domain $D$ which ends on the boundary of $D$. 
It can be thought of as interpreting the constraint that a 
SAW ends on the boundary of $D$ by thickening the boundary.
We weight a SAW $\omega$ by the usual weight of $\mu^{-|\omega|}$ 
times two factors. The factor $\weight(\omega)$, 
given by \reff{weight_def}, accounts for the varying thickness of 
the boundary. The other factor of $1/l(\theta)$ accounts for 
the local lattice effect near the endpoint of the SAW that 
persists in the scaling limit. 

We have conjectured that the scaling limit of the dilation ensemble
is related to the limit as $N \ra \infty$ of the uniform measure on 
SAW's with $N$-steps conditioned on the event that when the SAW
is scaled so that it ends on the boundary of $D$, the SAW lies entirely
in $D$. We have used our conjecture to simulate the SAW in several 
bounded domains and compared the distribution of the endpoint of the
SAW on the boundary with the SLE partition function prediction 
of this density. We find excellent agreement. 
This supports our conjectured relationship between 
the dilation ensemble and the uniform probability measure on SAW's of 
a fixed length,
the conjecture for the lattice effects that persist in the scaling
limit and the SLE partition function predictions for the boundary 
density. Note our simulations only  looked at the boundary density of 
the endpoint of the SAW. By using the SLE prediction for the boundary 
density and various explicit predictions about chordal and radial 
SLE, one can make explicit predictions about the SAW curve inside $D$
and use our conjecture to test them. 
We have not carried out such simulations.

\appendix

\section{Computation of SLE partition functions}

In this appendix we give the analytic computations of the boundary 
density for the various domains.
We use $\rho(\cdot)$ to denote the boundary density with respect 
to arc length (or a constant multiple of arc length) along $C$. 
In all cases the SAW starts at the origin. 
The constant $c$ in the following is determined by 
the requirement that $\rho$ is a probability density. 

\subsection{Strip - chordal}

The domain is a strip of height $1$ 
with the origin on the lower boundary. We are interested in the ensemble 
of SAW's in the strip that start at the origin and end at some point 
on the upper boundary of the strip. 
We parametrize the upper boundary by $x+i$ where $-\infty<x<\infty$. 
A straightforward application of \reff{conf_covar} yields the boundary 
density with respect to $dx$.
\beann
\rho(x) = c \left[\cosh(\pi x/2)\right]^{-5/4}
\eeann

\subsection{Strip - radial}

The domain is a horizontal strip of height $1$, with the 
origin at a height $h$ above the bottom boundary. So the domain is 
$\{z: -h \le Im(z) \le 1-h\}$. We consider the ensemble of SAW's 
in this strip that start at the origin, and end at a point 
on either the top or bottom boundary of the strip. 
We parametrize the lower boundary by $x-ih$ 
and the upper boundary by $x+(1-h)i$ where $-\infty<x<\infty$. 
Then the boundary density with respect to $dx$ is 
\beann
\rho(x+iy) = \cases{ 
       c \left[\cosh(\pi x) - \cos(\pi h)\right]^{-5/8}
         \quad if \quad y= -h \cr
       c \left[\cosh(\pi x) + \cos(\pi h)\right]^{-5/8}
         \quad if \quad y= 1-h \cr
}
\eeann

\subsection{Triangle}

We use the Schwarz-Christoffel mapping 
\beann
F(z) = \int_z^\infty \, (w+1)^{-2/3} \, (w-1)^{-2/3} \, dw 
\eeann
which maps the upper half plane onto an equilateral triangle. 
It sends $-1,1,\infty$ to the vertices, and some simple substitutions
show it sends $-3,0,3$ to the midpoints of the three sides. 
To find the pre-image of the center of the triangle, 
note that the center of triangle is the unique point fixed 
by rotations about the center. These rotations are the conformal 
automorphisms of the triangle that permute the vertices. So they 
correspond to conformal automorphisms of the half plane
that permute $-1,1$ and $\infty$. Such M\"obius transformations fix
$i \sqrt{3}$, so $F(i \sqrt{3})$ is the center of the triangle.

A simple application of \reff{conf_covar} shows that for SAW's in the 
half plane starting at $i \sqrt{3}$, the boundary density along the 
real axis is proportional to $(x^2+3)^{-5/8}$. 
By symmetry it is enough to find the boundary density for the triangle
between the midpoint of one edge and a vertex. 
We consider $F(3,\infty)$. (For $|z| \ge 3$, $F(z)$ can be computed
numerically very quickly by a power series expansion.)
Let $l_0$ be the length of the edges in 
the triangle. We define $l=l(x)$ by  
\beann
{2 l(x) \over l_0} = {F(x) \over F(3)}
\eeann
Then $l(x)$ is the distance from $F(x)$ to the vertex $F(\infty)$.
Applying \reff{conf_covar} we find the boundary density 
with respect to $dl$: 
\beann
\rho(l)= c {|F^{\prime}(x)|^{-5/8} \over (x^2+3)^{5/8}}=
c {(x^2-1)^{5/12} \over (x^2+3)^{5/8}}
\eeann

\subsection{Off-center circle}

The domain $D$ is a unit disc centered at $a+ib$ where  
$|a+ib|<1$ so that the origin is inside the disc. 
Let $\phi$ be the polar angle of a point on the boundary with respect to 
the center at $a+ib$. So $\phi$ is proportional to 
arc length along the boundary.
An application of \reff{conf_covar} with a Moibius transformation 
yields the boundary density with respect to $d \phi$.
\beann
\rho(\phi) = c \left[1+a^2+b^2 +2a\cos \phi +2b \sin \phi\right]^{-5/8}
\eeann

\subsection{Partial circle - chordal}

We consider the unit disc centered at $b$ where $b$ is real and 
$|b|<1$. So part of the disc lies below the real axis. 
We take the domain to be the intersection of this disc with 
the upper half plane.
We consider SAW's in this domain that start at the origin and 
end on the arc of the circle above the real axis. 

The circle intersects the real axis at $\pm d$ with $d = \sqrt{1-b^2}$.
The map $ -(z+d)/(z-d)$  maps the domain to the wedge 
$0 < \arg(z) < \beta$  where $\tan \beta = - \sqrt{1-b^2}/b$.
So if we let 
\beann
f(z) = \left[ - {z+d \over z-d} \right]^{\pi/\beta}
\eeann
then $f$ maps the domain to the upper half plane and 
sends $0$ to $1$. The endpoint of the walk is mapped to a point
on the negative real axis. 

An easy application of \reff{conf_covar} shows that for the SAW in 
the upper half plane starting at $1$ and ending at $x$ on the negative 
real axis, the boundary density is proportional to $(1-x)^{-5/4}$. 
So applying \reff{conf_covar} to the map $f$ shows that the 
boundary density with respect to $d \phi$ is 
\beann
\rho(\phi)= c \left[ {|z+d|^{-1+\pi/\beta} 
\over |z-d|^{1+\pi/\beta}} \right]^{5/8}
{1 \over (1-f(z))^{5/4}}
\eeann
where $z= ib+e^{i\phi}$.

\subsection{Tangent circle - chordal}

The domain is a disc of radius $1$ centered at $i$ so it is tangent 
to the real axis at the origin. 
Let $\phi$ be the polar angle with respect to $i$, so $\phi$ is proportional
to arc length along the boundary. 
We condition on the event that the walk stays in this 
disc and ends on the upper half of the circle bounding the disc, i.e., 
the arc of the boundary where $0 \le \phi \le 180$.

The boundary density for the upper half
plane when the SAW starts at the origin is $x^{-5/4}$. So 
using a M\"obius transformation 
to map the disc to the half plane, the formula
\reff{conf_covar} gives the density with respect to $d \phi$:
\beann
\rho(\phi) = c {|1-\cos \phi + \sin \phi|^{-5/4} (1-\cos \phi)^{5/8}}
\eeann

\end{document}